\documentclass[12pt]{article}
\usepackage{amssymb}
\newtheorem{theorem}{Theorem}[section]
\newtheorem{proposition}[theorem]{Proposition}
\newtheorem{lemma}[theorem]{Lemma}

\def\bull{\vrule height .9ex width .8ex depth -.1ex}
\newenvironment{proof}{\par\noindent {\bf Proof.~}}
              {\unskip\nobreak\hfill\hskip 2em \bull\par\medbreak}
\newenvironment{proofof}[1]{\smallbreak\noindent{\bf Proof of~#1.~}}
              {\unskip\nobreak\hfill\hskip 2em \bull\par\medbreak}
\def\om{\omega}
\def\eps{\varepsilon}
\def\La{\Lambda}
\def\la{\lambda}
\def\bR{\mathbb{R}}
\def\bC{\mathbb{C}}
\def\bZ{\mathbb{Z}}
\def\SL{\mathop{\mathrm{SL}}(2,\mathbb{R})}
\def\GL{\mathop{\mathrm{GL}}}
\def\cM{{\cal M}}
\def\tcM{\widetilde{\cal M}}
\def\cC{{\cal C}}
\def\cA{{\cal A}}
\def\cL{{\cal L}}
\def\cW{{\cal W}}
\def\tcW{\widetilde{\cal W}}
\def\cY{{\cal Y}}
\def\tcY{\widetilde{\cal Y}}
\def\cV{{\cal V}}
\def\cS{{\cal S}}
\def\Mod{\mathop{\mathrm{Mod}}}
\def\hol{\mathop{\mathrm{hol}}\nolimits}
\def\Aut{\mathop{\mathrm{Aut}}}
\title{Periodic geodesics on translation surfaces}
\author{Yaroslav Vorobets}
\date{July 15, 2003}
\begin{document}

\maketitle

\section{Introduction}\label{intro}

Let $M$ be a compact connected oriented surface.  The surface $M$ is called
a {\em translation surface\/} if it is equipped with a {\em translation
structure}, that is, an atlas of charts such that all transition functions
are translations in $\bR^2$.  It is assumed that the chart domains cover
all surface $M$ except for finitely many points called {\em singular}.  The
translation structure induces the structure of a smooth manifold, a flat
Riemannian metric, and a Borel measure on the surface $M$ punctured at the
singular points.  We require that the metric has a cone type singularity at
each singular point; then the area of the surface is finite.  The cone
angle is of the form $2\pi m$, where $m$ is an integer called the {\em
multiplicity\/} of the singular point.  A singular point of multiplicity
$1$ is called {\em removable\/}; it is rather a marked point than a true
singularity of the metric.

Furthermore, the translation structure allows us to identify the tangent
space at any nonsingular point $x\in M$ with the Euclidean space $\bR^2$.
In particular, the unit tangent space at any point is identified with the
unit circle $S^1=\{v\in\bR^2: |v|=1\}$.  The velocity is an integral of the
geodesic flow with respect to this identification.  Thus each oriented
geodesic has a {\em direction}, which is a uniquely determined vector in
$S^1$.  The direction of an unoriented geodesic is determined up to
multiplying by $\pm1$.

Suppose $X$ is a Riemann surface (one-dimensional complex manifold)
homeomorphic to the surface $M$.  Any nonzero Abelian differential on $X$
defines a translation structure on $M$.  The zeroes of the differential are
singular points of the translation structure, namely, a zero of order $k$
is a singular point of multiplicity $k+1$.  Every translation structure
without removable singular points can be obtained this way.

Any geodesic joining a nonsingular point to itself is {\em periodic\/} (or
{\em closed\/}).  We regard periodic geodesics as simple closed unoriented
curves.  Any periodic geodesic is included in a family of freely homotopic
periodic geodesics of the same length and direction.  The geodesics of the
family fill an open connected domain.  Unless the translation surface is a
torus without singular points, this domain is an annulus.  We call it a
{\em cylinder\/} of periodic geodesics (or simply a {\em periodic
cylinder\/}).  A periodic cylinder is bounded by geodesic segments of the
same direction with endpoints at singular points.  Such segments are called
{\em saddle connections}.

The fundamental results on periodic geodesics of translation surfaces were
obtained by Howard Masur in papers \cite{Masur86}, \cite{Masur88},
\cite{Masur90}.  These results can be summarized as follows.

\begin{theorem}[Masur]\label{Masur}
Let $M$ be a translation surface without removable singular points.

(a) There exists a periodic geodesic on $M$ of length at most $\alpha
\sqrt{S}$, where $S$ is the area of $M$ and $\alpha>0$ is a constant
depending only on the genus of $M$.

(b) The directions of periodic geodesics of $M$ are dense in $S^1$.

(c) Let $N_1(M,T)$ denote the number of periodic cylinders of $M$ of length
at most $T>0$.  Then there exist $0<c_1(M)<c_2(M)<\infty$ such that
$$
c_1(M)\le N_1(M,T)/T^2\le c_2(M)
$$
for $T$ sufficiently large.
\end{theorem}

The main goal of the present paper is to prove effective versions of
statements (a) and (c) of Theorem \ref{Masur}, and to generalize statement
(b).  In addition, we establish some properties of periodic geodesics on
generic translation surfaces.

Throughout the paper we consider translation surfaces that have at least
one singular point.  There is no loss of generality as we can declare an
arbitrary nonsingular point to be a removable singular point.

Our first result is an effective version of Theorem \ref{Masur}(a).

\begin{theorem}\label{main1}
Let $m$ be the sum of multiplicities of singular points of a translation
surface $M$, and $S$ be the area of $M$.  Then there exists a periodic
geodesic on $M$ of length at most $\alpha_m\sqrt{S}$, where $\alpha_m=
(8m)^{2^{3m-1}}$.
\end{theorem}

It should be admitted that the proof of Theorem \ref{Masur}(a) given by
Smillie in the survey \cite{S} can be further developed to obtain an
effective estimate of the constant $\alpha$ (unlike the proofs given in
\cite{Masur86} and \cite{MT}).  The techniques used below to prove Theorem
\ref{main1} are very similar to those used in \cite{S}.

The periodic geodesic provided by Theorem \ref{main1} belongs to a
cylinder of parallel periodic geodesics of the same length.  Although the
length of this cylinder is bounded, its width, in general, may be
arbitrarily small.  Nevertheless it is possible to find a periodic cylinder
whose area is not very small compared to the area of the whole surface.

\begin{theorem}\label{main2}
Let $m$ be the sum of multiplicities of singular points of a translation
surface $M$, and $S$ be the area of $M$.  Then there exists a cylinder of
periodic geodesics of length at most $\beta_m\sqrt{S}$, where $\beta_m=
2^{2^{4m}}$, and of area at least $S/m$.
\end{theorem}

The following theorem shows, in particular, that almost every point of a
translation surface lies on a periodic geodesic.

\begin{theorem}\label{bonus}
Let $m$ be the sum of multiplicities of singular points of a translation
surface $M$, and $S$ be the area of $M$.  For any $\delta\in(0,1)$ there
exist pairwise disjoint periodic cylinders $\La_1,\dots,\La_k$ of length at
most $(8m\delta^{-1})^{2^{3m-1}}\sqrt{S}$ such that the area of the union
$\La_1\cup\dots\cup\La_k$ is at least $(1-\delta)S$.
\end{theorem}

The group $\SL$ acts on the set of translation structures on a given
surface by postcomposition of the chart maps with linear transformations
from $\SL$.  This action preserves singular points along with their
multiplicities, geodesics, and the measure induced by translation
structure.  It does not preserve directions and lengths of geodesic
segments however.  This observation allows one to derive statement (b) of
Theorem \ref{Masur} from statement (a).  In the same way Theorem
\ref{main2} leads to the following result.

\begin{theorem}\label{main3}
Let $m$ be the sum of multiplicities of singular points of a translation
surface $M$, and $S$ be the area of $M$.  Then the directions of periodic
cylinders of area at least $S/m$ are dense in $S^1$.
\end{theorem}

A plane polygon is called {\em rational\/} if the angle between any two of
its sides is a rational multiple of $\pi$.  A construction of Zemlyakov and
Katok \cite{ZK} associates to any rational polygon $Q$ a translation
surface $M$ so that the study of the billiard flow in $Q$ can be reduced to
the study of the geodesic flow on $M$.  In view of this construction,
Theorem \ref{Masur}(b) implies that directions of periodic billiard orbits
in $Q$ are dense in the set of all directions.  Boshernitzan, Galperin,
Kr\"uger, and Troubetzkoy \cite{BGKT} strengthened this result.

\begin{theorem}[\cite{BGKT}]\label{4authors}
For any rational polygon $Q$, the periodic points of the billiard flow in
$Q$ are dense in the phase space of the flow.  Moreover, there exists a
dense $G_\delta$-set $Q_0\subset Q$ such that for every point $x\in Q_0$
the directions of periodic billiard orbits starting at $x$ form a dense
subset of $S^1$.
\end{theorem}

An analogous result for translation surfaces---periodic points of the
geodesic flow are dense in the phase space of the flow---can be obtained in
the same way (cf. \cite{MT}).  In this paper we prove a further
strengthening of Theorem \ref{Masur}(b).

\begin{theorem}\label{main4}
(a) For any translation surface $M$, there exists a $G_\delta$-set
$M_0\subset M$ of full measure such that for every point $x\in M_0$ the
directions of periodic geodesics passing through $x$ form a dense subset of
$S^1$.

(b) For any rational polygon $Q$, there exists a $G_\delta$-subset
$Q_0\subset Q$ of full measure such that for every point $x\in Q_0$ the
directions of periodic billiard orbits starting at $x$ form a dense subset
of $S^1$.
\end{theorem}

The next result is an effective version of Theorem \ref{Masur}(c).

\begin{theorem}\label{main5}
Let $M$ be a translation surface.  Denote by $N_1(M,T)$ the number of
cylinders of periodic geodesics on $M$ of length at most $T>0$.  By
$N_2(M,T)$ denote the sum of areas of these cylinders.  Then
$$
\Bigl((600m)^{(2m)^{2m}}\Bigr)^{-1}s^2S^{-2}T^2\le N_2(M,T)/S\le
N_1(M,T)\le (400m)^{(2m)^{2m}}s^{-2}T^2
$$
for any $T\ge 2^{2^{4m}}\sqrt{S}$, where $S$ is the area of $M$, $m$ is the
sum of multiplicities of singular points of $M$, and $s$ is the length of
the shortest saddle connection of $M$.
\end{theorem}

Let $M$ be a compact connected oriented surface of genus $p\ge1$.  For any
integer $n\ge1$ let $\cM(p,n)$ denote the set of equivalence classes of
isomorphic translation structures on $M$ with $n$ singular points (of
arbitrary multiplicity).  A point of $\cM(p,n)$ is a translation structure
considered up to isomorphism.  There is a natural topology on $\cM(p,n)$,
which is the topology of a locally compact metric space.  By $\cM_1(p,n)$
denote the subspace of $\cM(p,n)$ corresponding to translation structures
of area $1$.  The subspace $\cM_1(p,n)$ is endowed with a natural Borel
measure $\mu_0$, which is finite (see Section \ref{quad} for details).  In
general, the space $\cM_1(p,n)$ is not connected but the number of its
connected components is finite.  Let $\cC$ denote one of the connected
components.

For any translation structure $\om$ and any $T>0$ let $N_1(\om,T)$ denote
the number of periodic cylinders of $\om$ of length at most $T$.  By
$N_2(\om,T)$ denote the sum of areas of these cylinders.  The numbers
$N_1(\om,T)$ and $N_2(\om,T)$ do not change if we replace the translation
structure $\om$ by an isomorphic one.

\begin{theorem}\label{main6}
For $\mu_0$-a.e. $\om\in\cC$,
$$
\lim_{T\to\infty} N_1(\om,T)/T^2=c_1(\cC), \quad
\lim_{T\to\infty} N_2(\om,T)/T^2=c_2(\cC),
$$
where $c_1(\cC)$ and $c_2(\cC)$ are positive constants depending only on
the component $\cC$.
\end{theorem}

The first asymptotics in Theorem \ref{main6} was proved by Eskin and Masur
\cite{EM}.  The second asymptotics is obtained by applying results of
\cite{EM}.

The ratio $c_2(\cC)/c_1(\cC)$ may be regarded as the {\it mean area\/} of a
periodic cylinder of a generic area $1$ translation structure $\om\in\cC$.
It appears that $c_2(\cC)/c_1(\cC)=1/m_{\cC}$, where $m_{\cC}=2p-2+n$ is
the sum of multiplicities of singular points for translation structures in
$\cC$ (this will be proved in a subsequent paper).

Let $Y\to\cC$ be the fiber bundle over $\cC$ such that the fiber over a
point $\om\in\cC$ is the surface $M$ with the translation structure $\om$.
A point of $Y$ can be viewed as a pair $(\om,x)$, where $\om$ is a
representative of an equivalence class $\tilde\om\in\cC$ and $x\in M$ (the
point $x$ depends on the choice of $\om\in\tilde\om$) . The fiber bundle
$Y$ carries a natural finite measure $\mu_1$ that is the measure $\mu_0$
on the base $\cC$ and is the measure induced by translation structure on
the fiber.  Denote by $N_3(\om,x,T)$ the number of periodic geodesics of
a translation structure $\om$ of length at most $T$ that pass through a
point $x$.  This number does not change if we replace the pair $(\om,x)$ by
another representative of a point in $Y$.

\begin{theorem}\label{main7}
For $\mu_1$-a.e. $(\om,x)\in Y$,
$$
\lim_{T\to\infty} N_3(\om,x,T)/T^2=c_2(\cC),
$$
where the constant $c_2(\cC)$ is the same as in Theorem \ref{main6}.
\end{theorem}

The paper is organized as follows.  Section \ref{prel} contains
definitions, notation, and preliminaries.  The results on existence of
periodic geodesics (Theorems \ref{main1}, \ref{main2}, and \ref{bonus}) are
obtained in Section \ref{exist}.  The results on density of directions of
periodic geodesics (Theorems \ref{main3} and \ref{main4}) are obtained in
Section \ref{dense}.  Section \ref{lower} is devoted to the proof of
Theorem \ref{main5}.  In the final Section \ref{quad}, moduli spaces of
translation structures are considered.

\section{Preliminaries}\label{prel}

Let $M$ be a compact connected oriented surface.  A {\em translation
structure\/} on $M$ is an atlas of coordinate charts $\om=
\{(U_\alpha,f_\alpha)\}_{\alpha\in\cA}$, where $U_\alpha$ is a domain in
$M$ and $f_\alpha$ is a homeomorphism of $U_\alpha$ onto a domain in
$\bR^2$, such that:
\par$\bullet$ all transition functions are translations in $\bR^2$;
\par$\bullet$ chart domains $U_\alpha$, $\alpha\in\cA$, cover all surface
$M$ except for finitely many points (called {\em singular\/} points);
\par$\bullet$ the atlas $\om$ is maximal relative to the two preceding
conditions;
\par$\bullet$ a punctured neighborhood of any singular point covers a
punctured neighborhood of a point in $\bR^2$ via an $m$-to-$1$ map which is
a translation in coordinates of the atlas $\om$; the number $m$ is called
the {\em multiplicity\/} of the singular point.

A {\em translation surface\/} is a compact connected oriented surface
equipped with a translation structure.

The translation structures are also (and probably better) known as
``orientable flat structures'' or ``admissible positive $F$-structures''.

Let $M$ be a translation surface and $\om$ be the translation structure of
$M$.  Each translation of the plane $\bR^2$ is a smooth map preserving
orientation, Euclidean metric and Lebesgue measure on $\bR^2$.  Hence the
translation structure $\om$ induces a smooth structure, an orientation, a
flat Riemannian metric, and a finite Borel measure on the surface $M$
punctured at the singular points of $\om$.  Each singular point of $\om$ is
a cone type singularity of the metric.  The cone angle is equal to $2\pi
m$, where $m$ is the multiplicity of the singular point.  Any geodesic of
the metric is a straight line in coordinates of the atlas $\om$.  A
geodesic hitting a singular point is considered to be singular, its further
continuation is undefined.  Almost every element of the tangent bundle
gives rise to a nonsingular geodesic.

The translation structure $\om$ allows us to identify the tangent space at
any nonsingular point $x\in M$ with the Euclidean space $\bR^2$.  In
particular, the unit tangent space at any point is identified with the unit
circle $S^1=\{v\in\bR^2: |v|=1\}$.  The velocity is an integral of the
geodesic flow with respect to this identification.  Thus each oriented
geodesic is assigned a {\em direction\/} $v\in S^1$.  The direction of an
unoriented geodesic is determined up to multiplying by $\pm1$.  For any
$v\in S^1$, let $M_v$ denote the invariant surface of the phase space of
the geodesic flow corresponding to the movement with velocity $v$.
Clearly, the restriction of the geodesic flow to $M_v$ can be regarded as a
flow on the surface $M$.  This flow is called the {\em directional flow\/}
in direction $v$.  If a point $x\in M$ is singular or at least one of
geodesics starting at $x$ in the directions $\pm v$ hits a singular point,
then the directional flow is only partially defined at the point $x$.  The
directional flow is fully defined on a subset of full measure (depending on
$v$) and preserves the measure on $M$.

A singular point of multiplicity $1$ is called {\it removable}.  If $x\in
M$ is a removable singular point of the translation structure $\om$, then
there exists a translation structure $\om_+\supset\om$ such that $x$ is not
a singular point of $\om_+$.  On the other hand, let $x$ be a nonsingular
point of $M$.  Suppose $\om_-$ is the set of charts $(U_\alpha,f_\alpha)\in
\om$ such that $x\notin U_\alpha$.  Then $\om_-$ is a translation structure
on $M$ and $x$ is a removable singular point of $\om_-$.

Let $p$ be the genus of a translation surface $M$, $k$ be the number of
singular points of $M$, and $m$ be the sum of multiplicities of the
singular points.  Then $m=2p-2+k$.  It follows that there are no
translation structures on the sphere, a translation torus can have only
removable singular points, and a translation surface of genus $p>1$ has at
least one nonremovable singular point.

Suppose $X$ is a complex structure on a compact connected oriented surface
$M$.  Let $q$ be a nonzero Abelian differential (holomorphic 1-form) on
$X$.  A chart $(U,z)$, where $U$ is a domain in $M$ and $z$ is a
homeomorphism of $U$ onto a domain in $\bC$, is called a {\em natural
parameter\/} of the differential $q$ if $q=dz$ in $U$ with respect to the
complex structure $X$.  Let $\om$ denote the atlas of all natural
parameters of $q$.  The natural identification of $\bC$ with $\bR^2$ allows
us to consider $\om$ as an atlas of charts ranging in $\bR^2$.  It is easy
to observe that $\om$ is a translation structure on $M$.  The singular
points of $\om$ are the zeroes of the differential $q$, namely, a zero of
order $n$ is a singular point of multiplicity $n+1$.  Each translation
structure on $M$ without removable singular points can be obtained by this
construction.

Another way to construct translation surfaces is to glue them from
polygons.  Let $Q_1,\dots,Q_n$ be disjoint plane polygons.  The natural
orientation of $\bR^2$ induces an orientation of the boundary of every
polygon.  Suppose all sides of the polygons $Q_1,\dots,Q_n$ are groupped in
pairs such that two sides in each pair are of the same length and
direction, and of opposite orientations.  Glue the sides in each pair by
translation.  Then the union of the polygons $Q_1,\dots,Q_n$ becomes a
surface $M$.  By construction, the surface $M$ is compact and oriented.
Suppose $M$ is connected (if it is not, then we should apply the
construction to a smaller set of polygons).  The restrictions of the
identity map on $\bR^2$ to the interiors of the polygons $Q_1,\dots,Q_n$
can be regarded as charts of $M$.  This finite collection of charts extends
to a translation structure $\om$ on $M$.  The translation structure $\om$
is uniquely determined if we require that the set of singular points of
$\om$ be the set of points corresponding to vertices of the polygons
$Q_1,\dots,Q_n$.

A particular case of the latter construction is the so-called
Zemlyakov-Katok construction, which descends from the paper \cite{ZK}.  Let
$Q$ be a plane polygon.  Let $s_1,\dots,s_n$ be sides of $Q$.  For any $i$,
$1\le i\le n$, let $\tilde r_i$ denote the reflection of the plane in the
side $s_i$ and $r_i$ denote the linear part of $\tilde r_i$.  By $G(Q)$
denote the subgroup of $O(2)$ generated by the reflections $r_1,\dots,r_n$.
The polygon $Q$ is called {\em rational\/} if the group $G(Q)$ is finite.
All angles of a rational polygon are rational multiples of $\pi$.  This
property is equivalent to being rational provided the polygon is simply
connected.  Suppose the polygon $Q$ is rational.  Let $Q_g$, $g\in G(Q)$,
be disjoint polygons such that for any $g\in G(Q)$ there exists an isometry
$R_g:Q_g\to Q$ with linear part $g$.  Now for any $i\in\{1,\dots,n\}$ and
any $g\in G(Q)$ glue the side $R_g^{-1}s_i$ of the polygon $Q_g$ to the
side $R_{gr_i}^{-1}s_i$ of the polygon $Q_{gr_i}$ by translation.  This
transforms the union of polygons $Q_1,\dots,Q_n$ into a compact connected
oriented surface $M$.  Observe that the collection of isometries $R_g$,
$g\in G(Q)$, gives rise to a continuous map $f_Q:M\to Q$.  The surface $M$
is endowed with a translation structure $\om$ as described above.  Singular
points of $\om$ correspond to vertices of the polygon $Q$. Namely, the
vertex of any angle of the form $2\pi n_1/n_2$, where $n_1$ and $n_2$ are
coprime integers, gives rise to $N/n_2$ singular points of multiplicity
$n_1$, where $N$ is the cardinality of the group $G(Q)$.

The Zemlyakov-Katok construction is crucial for the study of the billiard
flow in rational polygons.  The {\em billiard flow\/} in a polygon $Q$ is a
dynamical system that describes a point-like mass moving freely within the
polygon $Q$ subject to elastic reflections in the boundary of $Q$.  A {\em
billiard orbit\/} in $Q$ is a broken line changing its direction at
interior points of the sides of $Q$ according to the law ``the angle of
incidence is equal to the angle of reflection''.  A billiard orbit hitting
a vertex of the polygon $Q$ is supposed to stop at this vertex.  A billiard
orbit starting at a point $x\in Q$ in a direction $v\in S^1$ is {\em
periodic\/} if it returns eventually to the point $x$ in the direction $v$.
Suppose $Q$ is a rational polygon.  Let $M$ be the translation surface
associated to $Q$.  Let $f_Q:M\to Q$ be the continuous map introduced
above.  It is easy to see that $f_Q$ maps any geodesic on the surface $M$
onto a billiard orbit in $Q$.  Conversely, any billiard orbit in $Q$ is the
image of a (not uniquely determined) geodesic on $M$.  By construction,
there exists a domain $D\subset M$ such that $f_Q$ maps the domain $D$
isometrically onto the interior of the polygon $Q$ and, moreover, the chart
$(D,f_Q|D)$ is an element of the translation structure of $M$.  If $L$ is a
geodesic starting at a point $x\in D$ in a direction $v\in S^1$, then
$f_Q(L)$ is the billiard orbit in $Q$ starting at the point $f_Q(x)$ in the
same direction.

Let $M$ be a translation surface.  A domain $U\subset M$ containing no
singular points is called a {\em triangle}\/ (a {\em polygon}\/, an
$n$-{\em gon}\/) if it is isometric to the interior of a triangle (resp. a
polygon, an $n$-gon) in the plane $\bR^2$.  Suppose $h:U\to P\subset\bR^2$
is a corresponding isometry.  The inverse map $h^{-1}:P\to U$ can be
extended to a continuous map of the closure of $P$ to $M$.  The images of
sides and vertices of the polygon $P$ under this extension are called {\em
sides}\/ and {\em vertices}\/ of $U$.  Every side of $U$ is either a
geodesic segment or a union of several parallel segments separated by
singular points.  Note that the number of vertices of the $n$-gon $U$ may
be less than $n$.  A {\it triangulation\/} of the translation surface $M$
is its partition into a finite number of triangles.

A {\em saddle connection\/} is a geodesic segment joining two singular
points or a singular point to itself and having no singular points in its
interior (note that singular points are saddles for directional flows).
Two saddle connections of a translation surface are said to be {\em
disjoint}\/ if they have no common interior points (common endpoints are
allowed).  Three saddle connections are pairwise disjoint whenever they are
sides of a triangle.  For any $T>0$ there are only finitely many saddle
connections of length at most $T$.  In particular, there exists the
shortest saddle connection (probably not unique).

The following proposition is well known (see, e.g. \cite{MT}, \cite{Vo}).

\begin{proposition}\label{prel1}
(a) Any collection of pairwise disjoint saddle connections can be extended
to a maximal collection.
\par (b) Any maximal collection of pairwise disjoint saddle connections
forms a triangulation of the surface $M$ such that all sides of each
triangle are saddle connections.
\par (c) For any maximal collection, the number of saddle connections is
equal to $3m$, and the number of triangles in the corresponding
triangulation is equal to $2m$, where $m$ is the sum of multiplicities of
singular points.
\end{proposition}

Any geodesic joining a nonsingular point to itself is called {\em
periodic\/} (or {\em closed\/}); such a geodesic is a periodic orbit of a
directional flow.  We only consider {\em primitive\/} periodic geodesics,
that is, the period of the geodesic is its length.  Also, we regard
periodic geodesics as unoriented curves.  If a geodesic starting at a point
$x\in M$ is periodic, then all geodesics starting at nearby points in the
same direction are also periodic.  Actually, each periodic geodesic belongs
to a family of freely homotopic periodic geodesics of the same length and
direction.  If $M$ is a torus without singular points, then this family
fills the whole surface $M$.  Otherwise the family fills a domain
homeomorphic to an annulus.  This domain is called a {\em cylinder\/} of
periodic geodesics (or simply a {\em periodic cylinder\/}) since it is
isometric to a cylinder $\bR/l\bZ\times(0,w)$, where $l,w>0$.  The numbers
$l$ and $w$ are called the {\em length\/} and the {\em width\/} of the
periodic cylinder.  The cylinder is bounded by saddle connections of the
same direction.

Let $\om=\{(U_\alpha,f_\alpha)\}_{\alpha\in\cA}$ be a translation structure
on the surface $M$.  For any linear operator $a\in\SL$ the atlas
$\{(U_\alpha,a{\circ}f_\alpha)\}_{\alpha\in\cal A}$ is also a translation
structure on $M$.  We denote this structure by $a\om$.  Clearly,
$(a_1a_2)\om=a_1(a_2\om)$ for any $a_1,a_2\in\SL$ so we have an action of
the group $\SL$ on the set of translation structures on $M$.  The
translation structures $\om$ and $a\om$ share the same singular points of
the same multiplicities and the same geodesics.  In addition, they induce
the same measure on the surface $M$.

To each oriented geodesic segment $L$ of the translation structure $\om$ we
associate the vector $v\in\bR^2$ of the same length and direction.  If the
segment $L$ is not oriented, then the vector $v$ is determined up to
reversing its direction.  For any $a\in\SL$ the vector $av$ is associated
to $L$ with respect to the translation structure $a\om$.  Given a direction
$v_1\in S^1$, the length of the orthogonal projection of $v$ on the
direction $v_1$ is called the {\em projection\/} of the segment $L$ on
$v_1$ (with respect to $\om$).

\section{Existence of periodic geodesics}\label{exist}

To prove Theorems \ref{main1}, \ref{main2}, and \ref{bonus}, we need the
following proposition.

\begin{proposition}\label{exist1}
Let $M$ be a translation surface of area $S$.  Suppose $L_1,\dots,L_k$
$(k\ge0)$ are pairwise disjoint saddle connections of length at most
$\sqrt{2S}$.  Then at least one of the following possibilities occur:
\par (1) saddle connections $L_1,\ldots,L_k$ partition the surface into
finitely many domains such that each domain is either a periodic cylinder
of length at most $\sqrt{S}$ or a triangle bounded by three saddle
connections;
\par (2) there exists a saddle connection $L$ of length at most
$2\sqrt{2S}$ disjoint from $L_1,\ldots,L_k$.
\end{proposition}

\begin{proof}
First consider the case when a small neighborhood of some singular point
$x_0$ contains an open sector $K$ of angle $\pi$ that is disjoint from
saddle connections $L_1,\ldots,L_k$.  Suppose there exists a geodesic
segment $J$ of length at most $\sqrt{2S}$ that goes out of the point $x_0$
across sector $K$ and ends in a point $y$ which is either a singular point
or an interior point of some saddle connection $L_j$, $1\le j\le k$.  We
can assume without loss of generality that the interior of the segment $J$
contains no singular point and is disjoint from saddle connections
$L_1,\ldots,L_k$.  If $y$ is a singular point, then $J$ is a saddle
connection disjoint from $L_1,\dots,L_k$, thus condition (2) holds.
Suppose $y$ is an interior point of $L_j$.  Let $K'$ be an open sector with
vertex at the point $x_0$ crossed by the segment $J$.  We assume that each
geodesic $I$ going out of $x_0$ across the sector $K'$ intersects $L_j$
before this geodesic hits a singular point or intersects another given
saddle connection.  This condition holds, for instance, when the angle of
the sector $K'$ is small enough.  Let $I_0$ denote the segment of the
geodesic $I$ from the point $x_0$ to the first intersection with $L_j$.
The segments $I_0$, $J$, and a subsegment of $L_j$ are sides of a triangle,
hence the length of $I_0$ is less than the sum of lengths of $J$ and $L_j$,
which, in turn, is at most $2\sqrt{2S}$.  Without loss of generality it can
be assumed that $K'$ is the maximal sector with the above property.  By the
maximality, both geodesics going out of $x_0$ along the boundary of $K'$
hit singular points before they intersect any of the given saddle
connections.  It follows that these geodesics are saddle connections of
length at most $2\sqrt{2S}$.  By construction, any of the two saddle
connections either is disjoint from saddle connections $L_1,\ldots,L_k$ or
coincides with one of them.  Since the angle of the sector $K'$ is less
than $\pi$, at least one of the boundary saddle connections crosses the
sector $K$; such a saddle connection is not among $L_1,\dots,L_k$.  Thus
condition (2) holds.

Now suppose that any geodesic segment of length $\sqrt{2S}$ going out of
the point $x_0$ across sector $K$ does not reach a singular point and does
not intersect saddle connections $L_1,\dots,L_k$.  Let $I$ be the geodesic
segment of length $\sqrt{S}$ that goes out of the singular point $x_0$
dividing the sector $K$ into two equal parts.  By $v$ denote one of two
directions orthogonal to the direction of $I$.  Let $\{F^t\}_{t\in\bR}$ be
the directional flow in direction $v$.  By $I_0$ denote the segment $I$
without its endpoints.  For any $t>0$ let $D_t$ denote the set of points of
the form $F^\tau x$, where $x\in I_0$ and $0<\tau<t$.  If $t$ is small
enough, then $D_t$ is a rectangle with sides $\sqrt{S}$ and $t$.  Let $t_1$
be the maximal number with this property.  The area of the rectangle
$D_{t_1}$ is equal to $t_1\sqrt{S}$, hence $t_1\le\sqrt{S}$.  Set
$D'_{t_1}=D_{t_1}\cup F^{t_1}I_0$.  Any point $x\in D'_{t_1}$ can be joined
to the point $x_0$ by a geodesic segment $J_x$ such that all interior
points of $J_x$ are contained in $D_{t_1}$.  The segment $J_x$ crosses the
sector $K$ and the length of $J_x$ is at most $\sqrt{t_1^2+S}\le\sqrt{2S}$.
It follows that the set $D'_{t_1}$ is disjoint from saddle connections
$L_1,\dots,L_k$ and contains no singular points.  Likewise, the set
$D'_{-t_1}=\{F^tx\mid x\in I_0,\ -t_1\le t<0\}$ is also disjoint from
$L_1,\dots,L_k$ and contains no singular points.  By the choice of $t_1$,
there exists a point $x_1\in I_0$ such that the point $F^{t_1}x_1$ either
is singular or belongs to $I_0$.  As $F^{t_1}x_1\in D'_{t_1}$, we have
$F^{t_1}x_1\in I_0$.  Let $y_0$ denote the endpoint of $I$ different from
$x_0$.  It is easy to see that the set $I_1=I_0\cap F^{t_1}I_0$ is an open
subsegment of $I_0$ and $y_0$ is an endpoint of $I_1$.  Suppose that
$I_1\ne I_0$.  Let $y$ be an endpoint of $I_1$ that is an interior point of
$I$.  Obviously, $F^{-t_1}y\in I$.  Since $y$ is an endpoint of $I_1$, it
follows that $F^{-t_1}y=x_0$.  On the other hand, $F^{-t_1}y\in D'_{-t_1}$.
This contradiction proves that $I_1=I_0$.  It follows that the set
$D'_{t_1}$ is a union of periodic geodesics of length $t_1$ and of
direction $v$.  Therefore $D'_{t_1}$ is contained in a periodic cylinder
$\La$ of length $t_1\le\sqrt{S}$.  The cylinder $\La$ contains the sector
$K$.  By construction, at least some of periodic geodesics in $\La$ do not
intersect saddle connections $L_1,\ldots,L_k$.  It follows easily that the
whole cylinder is disjoint from $L_1,\dots,L_k$.

Now suppose that condition (2) does not hold.  We have to show that
condition (1) does hold in this case.  By the above the saddle connections
$L_1,\dots,L_k$ divide a small neighborhood of each singular point into
sectors of angle at most $\pi$.  Moreover, each sector of angle $\pi$ is
contained within a periodic cylinder of length at most $\sqrt{S}$ disjoint
from $L_1,\dots,L_k$.  The saddle connections $L_1,\dots,L_k$ partition the
surface $M$ into finitely many domains.  Let $D$ be one of these domains.
Take a singular point $x_0$ at the boundary of $D$.  A small neighborhood
of $x_0$ intersects the domain $D$ in one or more sectors of angle at most
$\pi$.  Let $K$ be one of such sectors.  If the sector $K$ is of angle
$\pi$, then it is contained in a periodic cylinder $\La$ of length at most
$\sqrt{S}$ disjoint from $L_1,\dots,L_k$.  Clearly, $\La\subset D$.  The
lengths of saddle connections bounding the cylinder $\La$ do not exceed the
length of $\La$.  Any of these saddle connections either is disjoint from
$L_1,\dots,L_k$ or coincides with one of them.  Since condition (2) does
not hold, all saddle connections bounding $\La$ are among $L_1,\dots,L_k$.
This means that $D=\La$.  Now consider the case when the angle of the
sector $K$ is less than $\pi$.  The sector $K$ is bounded by some saddle
connections $L_i$ and $L_j$.  Let $T$ be a triangle such that $T$ contains
the sector $K$, the saddle connection $L_i$ is a side of $T$, and a
subsegment of $L_j$ is another side of $T$.  Obviously, $T\subset D$.  We
can assume without loss of generality that $T$ is the maximal triangle with
this property.  By $L_0$ denote the side of $T$ different from $L_i$ and
from the subsegment of $L_j$.  Let $J$ be a geodesic segment that goes out
of the point $x_0$ across sector $K$, crosses the triangle $T$, and ends in
a point $y\in L_0$.  The length of $J$ is less than the sum of lengths of
$L_i$ and $L_j$, which, in turn, is at most $2\sqrt{2S}$.  Since condition
(2) does not hold, the point $y$ can not be singular.  It follows that the
side $L_0$ is a single geodesic segment.  By the maximality of $T$, the
whole saddle connection $L_j$ is a side of $T$.  Then $L_0$ is a saddle
connection.  By the triangle inequality, the length of $L_0$ is at most
$2\sqrt{2S}$.  Hence $L_0$ is one of the saddle connections $L_1,\dots,L_k$
as otherwise $L_0$ is disjoint from $L_1,\dots,L_k$.  This means that
$D=T$.  Thus condition (1) holds.
\end{proof}

For any operator $a\in\SL$, let $\|a\|$ denote the {\em norm\/} of $a$
and $C(a)$ denote the {\em condition number}\/ of $a$:
$$
\|a\|=\max\limits_{v\in\bR^2,\,|v|=1}|av|, \qquad
C(a)=\max\,(\|a\|,\|a^{-1}\|).
$$
Obviously, $C(a_1a_2)\le C(a_1)C(a_2)$ for any $a_1,a_2\in\SL$.  Suppose
$L$ is a geodesic segment of a translation structure $\om$.  Then the
lengths of the segment $L$ with respect to translation structures $\om$ and
$a\om$ differ by at most $C(a)$ times.

\begin{proofof}{Theorem \ref{bonus}}
Let $\delta\in(0,1)$.  Set $\eps=(8m\delta^{-1})^{-2^{3m-1}}$.  Note that
$\eps^{2^{-3m}}\le1/2$.  Let $\om$ denote the translation structure of the
translation surface $M$.  Suppose that a sequence $L_1,\dots,L_k$ ($k\ge0$)
of pairwise disjoint saddle connections and a sequence of operators
$a_0,a_1,\dots,a_k\in\SL$ satisfy the following two conditions: (i)
$C(a_i)\le (1/\eps)^{1-2^{-i}}$ for $i=0,1,\dots,k$; and (ii) the length of
$L_i$ with respect to the translation structure $a_j\om$, $1\le i\le j\le
k$, does not exceed $2\sqrt{2S}\,\eps^{2^{-j}}$.  Furthermore, suppose
there exists a saddle connection $L$ disjoint from $L_1,\dots,L_k$ and of
length at most $2\sqrt{2S}$ with respect to the translation structure
$a_k\om$.  Let $v,u\in S^1$ be orthogonal vectors such that $v$ is parallel
to the saddle connection $L$ with respect to $a_k\om$.  Define an operator
$b\in\SL$ by equalities $bv=\eps^{2^{-k-1}}v$, $bu=\eps^{-2^{-k-1}}u$.
Further, set $a_{k+1}=ba_k$.  Obviously, $C(b)=(1/\eps)^{2^{-k-1}}$, hence
$C(a_{k+1})\le C(b)C(a_k)\le (1/\eps)^{2^{-k-1}}(1/\eps)^{1-2^{-k}}=
(1/\eps)^{1-2^{-k-1}}$.  The length of the saddle connection $L$ with
respect to $a_{k+1}\om$ is at most $2\sqrt{2S}\,\eps^{2^{-k-1}}$, while the
length of saddle connections $L_1,\ldots,L_k$ with respect to $a_{k+1}\om$
is at most $2\sqrt{2S}\,\eps^{2^{-k}}C(b)=2\sqrt{2S}\,\eps^{2^{-k-1}}$.
Thus the sequence of saddle connections $L_1,\dots,L_k,L$ and the sequence
of operators $a_0,a_1,\dots,a_k,a_{k+1}$ satisfy the conditions (i) and
(ii).

Pairs of sequences satisfying conditions (i) and (ii) do exist, for
example, the empty sequence of saddle connections and the sequence
consisting of one operator $a_0=1$.  By Proposition \ref{prel1}, the number
of pairwise disjoint saddle connections can not exceed $3m$.  Therefore
there exists a pair of sequences $L_1,\dots,L_k$ and $a_0,a_1,\dots,a_k$
satisfying conditions (i) and (ii) with maximal possible $k$.  The lengths
of the saddle connections $L_1,\dots,L_k$ with respect to $a_k\om$ are at
most $2\sqrt{2S}\eps^{2^{-k}}\le 2\sqrt{2S}\eps^{2^{-3m}}\le \sqrt{2S}$,
thus Proposition \ref{exist1} applies.  By the maximality of $k$, there is
no saddle connection disjoint from $L_1,\ldots,L_k$ and of length at most
$2\sqrt{2S}$ with respect to the translation structure $a_k\om$.  Thus the
saddle connections $L_1,\ldots,L_k$ partition the surface $M$ into finitely
many domains such that each domain is either a periodic cylinder of length
at most $\sqrt{S}$ with respect to $a_k\om$ or a triangle bounded by three
saddle connections.  Any triangle in this partition is of area at most
$\frac12(2\sqrt{2S}\,\eps^{2^{-k}})^2$ with respect to both $a_k\om$ and
$\om$.  It follows from Proposition \ref{prel1} that there are at most $2m$
triangles in the partition.  Hence the union of these triangles is of area
at most
$$
m(2\sqrt{2S}\,\eps^{2^{-k}})^2\le m(2\sqrt{2S}\,\eps^{2^{-3m}})^2=\delta S.
$$
Then the union of all periodic cylinders in the partition is of area at
least $(1-\delta)S$.  It remains to observe that the length of each
periodic cylinder with respect to the translation structure $\om$ is at
most $C(a_k)\sqrt{S}\le (1/\eps)^{1-2^{-k}}\sqrt{S}< \eps^{-1}\sqrt{S}=
(8m\delta^{-1})^{2^{3m-1}}\sqrt{S}$.
\end{proofof}

\begin{proofof}{Theorem \ref{main1}}
By Theorem \ref{bonus}, for any $\delta\in(0,1)$ the translation surface
$M$ admits a periodic geodesic of length at most
$(8m\delta^{-1})^{2^{3m-1}}\sqrt{S}$.  For any $T>0$ the number of periodic
cylinders of length at most $T$ is finite, therefore there exists a
shortest periodic geodesic.  Let $l$ denote its length.  Since $l\le
(8m\delta^{-1})^{2^{3m-1}}\sqrt{S}$ for any $\delta\in(0,1)$, we have $l\le
(8m)^{2^{3m-1}}\sqrt{S}$.
\end{proofof}

\begin{proofof}{Theorem \ref{main2}}
In the case $m=1$ the translation surface $M$ is a torus with one removable
singular point.  Here every cylinder of periodic geodesics fills the whole
surface (up to the boundary saddle connection).  As $\alpha_1=8^{2^2}=
2^{12}<2^{2^4}=\beta_1$, the theorem follows from Theorem \ref{main1} in
this case.

Consider the case $m\ge2$.  By Theorem \ref{bonus}, there exist pairwise
disjoint periodic cylinders $\La_1,\dots,\La_k$ of length at most
$(8m^2)^{2^{3m-1}}\sqrt{S}$ such that the area of the union
$\La_1\cup\dots\cup\La_k$ is at least $(1-1/m)S$.  Each cylinder $\La_i$
can be triangulated by pairwise disjoint saddle connections.  The number of
triangles in any triangulation is at least $2$.  If we require that each
side of any triangle is a saddle connection (not a union of several saddle
connections), then the number of triangles is equal to the number of saddle
connections bounding $\La_i$, where saddle connections bounding $\La_i$
from both sides should be counted twice.  All saddle connections used in
triangulation of the cylinders $\La_1,\dots,\La_k$ are pairwise disjoint
since the cylinders are disjoint.  By Proposition \ref{prel1}, we can add
several saddle connections to obtain a partition of the surface $M$ into
$2m$ triangles bounded by disjoint saddle connections.  It follows easily
that the number $k$ of cylinders is at most $m$.  Moreover, if $k=m$ then
the closure of the union $\La_1\cup\dots\cup\La_k$ is the whole surface
$M$.  In the latter case at least one of the cylinders $\La_1,\dots,\La_k$
is of area not less than $S/m$.  In the case $k<m$ one of the cylinders is
of area not less than $(1-1/m)S/k\ge (1-1/m)S/(m-1)=S/m$.

To complete the proof, it remains to show that $(8m^2)^{2^{3m-1}}\le
\beta_m$.  It is easy to observe that $8m^2\le 2^{2m+1}$ and $2m+1<2^{m+1}$
for any integer $m\ge1$.  Hence, $(8m^2)^{2^{3m-1}}\le 2^{(2m+1)2^{3m-1}}<
2^{2^{4m}}=\beta_m$.
\end{proofof}

\section{Density of directions of periodic geodesics}\label{dense}

In this section we prove Theorems \ref{main3} and \ref{main4}.  They are
derived from Theorems \ref{main2} and \ref{bonus}, respectively.

\begin{proofof}{Theorem \ref{main3}}
Let $\om$ denote the translation structure of the translation surface $M$.
We have to show that for any direction $v\in S^1$ and any $\eps>0$ there
exists a periodic cylinder of $\om$ of area at least $S/m$ such that the
angle between $v$ and the direction of the cylinder is less than $\eps$.

Let $u\in S^1$ be a direction orthogonal to $v$.  For any $\la>1$ define an
operator $a_\la\in\SL$ by equalities $a_\la v=\la^{-1}v$, $a_\la u=\la u$.
The sum of multiplicities of singular points of the translation structure
$a_\la\om$ is equal to $m$ and the area of the surface $M$ with respect to
$a_\la\om$ is equal to $S$.  By Theorem \ref{main2}, there exists
a periodic cylinder $\La_\la$ of area at least $S/m$ such that the length
of $\La_\la$ with respect to $a_\la\om$ is at most $l=2^{2^{4m}}\sqrt{S}$.
Let $h_\la$ and $w_\la$ be projections of a periodic geodesic from the
cylinder $\La_\la$ on the directions $v$ and $u$ (with respect to the
translation structure $\om$).  Further, let $\varphi_\la$ be the angle
between $v$ and the direction of $\La_\la$, where the direction of the
cylinder is chosen so that $0\le\varphi_\la\le\pi/2$.  Obviously,
$h_\la\le\la l$, $w_\la\le\la^{-1}l$.  Let $s$ denote the length of the
shortest saddle connection of $\om$.  Assuming $\la$ is large enough, we
have $w_\la\le s/\sqrt{2}$.  Since the length of the cylinder $\La_\la$,
which is equal to $\sqrt{h^2_\la+w^2_\la}$, is not less than $s$, it
follows that $h_\la\ge s/\sqrt{2}$.  Then $\varphi_\la\le \tan\varphi_\la=
w_\la/h_\la\le \la^{-1}l\sqrt{2}/s$, which tends to zero as $\la$ goes to
infinity.
\end{proofof}

\begin{proofof}{Theorem \ref{main4}}
Let $M$ be a translation surface.  Take a nonempty open subset $U$ of the
circle $S^1$.  Let $P_U$ denote the set of points $x\in M$ lying on
periodic geodesics with directions in the set $U$.  The set $P_U$ is open.
Let us show that this set is of full measure.  Take a vector $v\in U$.
Choose $\eps>0$ such that a direction $v'\in S^1$ belongs to $U$ whenever
the angle between $v'$ and $v$ is less than $\eps$.  Further, choose some
$\delta\in(0,1)$.  Let $u\in S^1$ be a direction orthogonal to $v$.  For
any $\la>1$ define an operator $a_\la\in\SL$ by equalities $a_\la v=
\la^{-1}v$, $a_\la u=\la u$.  Let $\om$ denote the translation structure of
$M$, $m$ denote the sum of multiplicities of singular points of $M$, and
$S$ denote the area of $M$.  By Theorem \ref{bonus}, there exist pairwise
disjoint periodic cylinders $\La_1,\dots,\La_k$ such that the length of
every cylinder with respect to the translation structure $a_\la\om$ is at
most $l_\delta=(8m\delta^{-1})^{2^{3m-1}}\sqrt{S}$ and the area of the
union $\La_1\cup\dots\cup\La_k$ is at least $(1-\delta)S$ (with respect to
both $a_\la\om$ and $\om$).  Take some cylinder $\La_i$.  Let $h$ and $w$
be projections of a periodic geodesic from the cylinder $\La_i$ on the
directions $v$ and $u$ (with respect to the translation structure $\om$).
Further, let $\varphi$ be the angle between $v$ and the direction of
$\La_i$ ($0\le\varphi\le\pi/2$).  Obviously, $h\le\la l_\delta$,
$w\le\la^{-1}l_\delta$.  Let $s$ denote the length of the shortest saddle
connection of $\om$.  The length of the cylinder $\La_i$ is not less than
$s$.  If $\la\ge l_\delta\sqrt{2}/s$, then $w\le s/\sqrt{2}$, hence $h\ge
s/\sqrt{2}$.  It follows that $\varphi\le \tan\varphi=w/h\le \la^{-1}
l_\delta\sqrt{2}/s$.  If, moreover, $\la>\eps^{-1}l_\delta\sqrt{2}/s$, then
$\varphi<\eps$ and the direction of the cylinder $\La_i$ is in the set $U$.
Thus the cylinders $\La_1,\dots,\La_k$ are contained in the set $P_U$
provided $\la$ is sufficiently large.  Then the area of $P_U$ is at least
$(1-\delta)S$.  As $\delta$ can be chosen arbitrarily small, the area of
$P_U$ is equal to $S$.

Choose a sequence $U_1,U_2,\dots$ of nonempty open subsets of the circle
$S^1$ such that any other nonempty open subset of $S^1$ contains some
$U_i$.  By the above the sets $P_{U_1},P_{U_2},\dots$ are open sets of full
measure.  Hence the set $P_\infty=\cap_{i=1}^\infty P_{U_i}$ is a
$G_\delta$-subset of full measure of the surface $M$.  Take a point $x\in
P_\infty$.  For any positive integer $i$ there exists a direction $v_i\in
U_i$ that is the direction of a periodic geodesic passing through $x$.  By
construction, the sequence $v_1,v_2,\dots$ is dense in $S^1$.  The first
statement of the theorem is proved.

To derive statement (b) of Theorem \ref{main4} from statement (a), we only
need to recall the Zemlyakov-Katok construction (see Section \ref{prel}).
Let $Q$ be a rational polygon and $M$ be the translation surface associated
to $Q$.  By construction, there is a continuous map $f:M\to Q$ and a domain
$D\subset M$ containing no singular points such that $f$ maps the domain
$D$ isometrically onto the interior of the polygon $Q$.  Moreover, if $L$
is a geodesic passing through a point $x\in D$ in a direction $v\in S^1$,
then $f(L)$ is the billiard orbit in $Q$ starting at the point $f(x)$ in
the direction $v$.  The billiard orbit $f(L)$ is periodic if and only if
the geodesic $L$ is periodic.  By the above there exists a $G_\delta$-set
$M_0\subset M$ of full measure such that for any $x\in M_0$ the directions
of periodic geodesics passing through $x$ are dense in $S^1$.  Then the set
$Q_0=f(D\cap M_0)$ is a $G_\delta$-subset of the polygon $Q$ and the area
of $Q_0$ is equal to the area of $Q$.  For any point $x\in Q_0$ the
directions of periodic billiard orbits in $Q$ starting at $x$ are dense in
$S^1$.
\end{proofof}

\section{Lower quadratic estimates}\label{lower}

Let $M$ be a translation surface.  In this section we obtain effective
estimates of the growth functions $N_1(M,\cdot)$ and $N_2(M,\cdot)$, where
$N_1(M,T)$ is the number of cylinders of periodic geodesics of length at
most $T$ and $N_2(M,T)$ is the sum of areas of those cylinders.  Throughout
the section $m$ denotes the sum of multiplicities of singular points of the
translation surface $M$ ($m\ge1$), $S$ denotes the area of $M$, and $s$
denotes the length of the shortest saddle connection of $M$.

For any $T>0$ let $N_0(M,T)$ denote the number of saddle connection of $M$
of length at most $T$.  An effective upper estimate of this number was
obtained in \cite{Vo}.

\begin{theorem}[\cite{Vo}]\label{lower1}
$N_0(M,T)\le h_ms^{-2}T^2$ for any $T>0$, where $h_1=(3\cdot 2^7)^6$ and
$h_m=(400m)^{(2m)^{2m}}$ for $m\ge2$.
\end{theorem}

In the case $m=1$ the latter estimate can be significantly improved.

\begin{lemma}\label{lower2}
Suppose $M$ is a translation torus with a single (removable) singular
point.  Then $N_0(M,T)\le 7s^{-2}T^2$ for any $T>0$.  In addition, $s^2\le
3S/2$.
\end{lemma}

\begin{proof}
The translation torus $M$ is isometric to a torus $\bR^2/(v_1\bZ\oplus
v_2\bZ)$, where $v_1$ and $v_2$ are linearly independent vectors in $\bR^2$
(we do not require that the isometry preserve directions).  Let $\cL=
v_1\bZ\oplus v_2\bZ$.  By $\cL_0$ denote the set of vectors in $\cL$
contained in neither of the lattices $2\cL,3\cL,\dots$.  The isometry
establishes a one-to-one correspondence between saddle connections of $M$
and pairs of vectors $\pm v\in\cL_0$.  The length of a saddle connection is
equal to the length of the corresponding vectors.  By $H$ denote the set of
points $(y_1,y_2)\in\bR^2$ such that either $y_2>0$, or $y_2=0$ and
$y_1>0$.  Then $N_0(M,T)$ is equal to the number of vectors of length at
most $T$ in the set $H\cap\cL_0$.  Notice that the vectors $v_1$ and $v_2$
are not determined in a unique way.  Without loss of generality it can be
assumed that $v_1=(0,s)$ and $v_2=(S/s,y)$, where $0\le y<s$.  Then
$\min(|v_2|,|v_2-v_1|)\le\sqrt{(S/s)^2+(s/2)^2}$.  Since $|v|\ge s$ for any
nonzero vector $v\in v_1\bZ\oplus v_2\bZ$, we have $s^2\le(S/s)^2+(s/2)^2$.
It follows that $s^2\le 2S/\sqrt{3}\le 3S/2$.

Let $i$ and $j$ be positive integers.  The rectangle $P^+_{i,j}=
((i-1)S/s,iS/s]\times[(j-1)s,js)\subset\bR^2$ is contained in the halfplane
$H$ and contains precisely one element of the lattice $\cL$.  Likewise, the
rectangle $P^-_{i,j}=[-iS/s,-(i-1)S/s)\times((j-1)s,js]$ is contained in
$H$ and contains only one element of $\cL$.  Given $T>0$, let $B_T=
\{v\in\bR^2: |v|\le T\}$.  The number of rectangles of the form
$P^\pm_{i,j}$ contained in the halfdisc $B_T\cap H$ does not exceed $\pi
T^2/(2S)$.  For any $k\in\bZ$ set $L_k=\{(y_1,y_2)\in H: y_1=Sk/s\}$.  If
$k\ne0$, then the halfline $L_k$ containes at most one element $v\in\cL$
such that $v\in B_T$ but the rectangle of the form $P^\pm_{i,j}$ containing
$v$ is not contained in $B_T$.  The halfline $L_0$ containes at most $T/s$
elements of $B_T\cap\cL$.  Finally, the cardinality of the set $B_T\cap
H\cap\cL$ is at most $\pi T^2/(2S)+2sT/S+T/s$.  As this cardinality is not
less than $N_0(M,T)$, we have $N_0(M,T)\le \pi T^2/(2S)+2sT/S+T/s\le 3\pi
s^{-2}T^2/4+4s^{-1}T\le 3s^{-2}T^2+4s^{-1}T$.  It follows that $N_0(M,T)\le
7s^{-2}T^2$ for $T\ge s$.  If $T<s$, then $N_0(M,T)=0<7s^{-2}T^2$.
\end{proof}

The following lemma is an improved version of Theorem \ref{main2} for
translation tori.

\begin{lemma}\label{lower3}
Suppose $M$ is a translation torus with $m\ge1$ singular points.  Then
there exists a cylinder of periodic geodesics of length at most $2\sqrt{S}$
and of area at least $S/m$.
\end{lemma}

\begin{proof}
Let $\om$ denote the translation structure of the translation torus $M$.
Let $x_1,x_2,\dots,x_m$ be the singular points of $\om$.  All singular
points are removable.  By $\om_1$ denote the translation structure on $M$
such that $\om_1\supset\om$ and $x_1$ is the only singular point of
$\om_1$.  Let $M_1$ denote the torus $M$ considered as the translation
surface with the translation structure $\om_1$.  Let $S_1$ be the area of
$M_1$ and $s_1$ be the length of the shortest saddle connection of $M_1$.
It is easy to observe that $S_1=S$ and $s_1\ge s$.  The shortest saddle
connection of $M_1$ bounds a periodic cylinder $\La$ of $M_1$.  The length
of $\La$ is equal to $s_1$ and the area of $\La$ is equal to $S$.  By Lemma
\ref{lower2}, $s_1\le \sqrt{3S/2}<2\sqrt{S}$.  The points $x_2,\dots,x_m$
split the cylinder $\La$ into several periodic cylinders of the translation
surface $M$.  All these cylinders are of length $s_1\le 2\sqrt{S}$.  The
number of the cylinders does not exceed $m$, hence at least one of them is
of area not less than $S/m$.
\end{proof}

\begin{proofof}{Theorem \ref{main5}}
Let $M$ be a translation surface.  To each cylinder of periodic geodesics
of $M$ we assign a saddle connection bounding the cylinder.  The length of
the saddle connection does not exceed the length of the cylinder.  It is
possible that a saddle connection bounds two different periodic cylinders.
Nevertheless the assignment can be done so that any saddle connection is
assigned to at most one cylinder.  It follows that $N_1(M,T)\le N_0(M,T)$
for any $T>0$.  Thus Theorem \ref{lower1} (in the case $m\ge2$) and Lemma
\ref{lower2} (in the case $m=1$) imply that $N_1(M,T)\le (400m)^{(2m)^{2m}}
s^{-2}T^2$ for any $T>0$.  Besides, the estimate $N_2(M,T)/S\le N_1(M,T)$
is trivial.

We proceed to the proof of the lower estimate of $N_2(M,T)$.  Let
$C_0=\tilde h_m s^{-2}$, where $\tilde h_1=7$ and $\tilde h_m=
(400m)^{(2m)^{2m}}$ for $m\ge2$.  By Theorem \ref{lower1} and Lemma
\ref{lower2}, $N_1(M,T)\le N_0(M,T)\le C_0T^2$ for any $T>0$.  Denote by
$\sigma(T)$ the sum of inverse lengths over all cylinders of periodic
geodesics of length at most $T$.  Let $T_1\le T_2\le\ldots\le T_n\le\ldots$
be lengths of periodic cylinders of $M$ in ascending order.  It follows
from the estimate $N_1(M,T)\le C_0T^2$ that $T_n\ge C_0^{-1/2}n^{1/2}$,
$n=1,2,\ldots$.  Therefore,
$$
\sigma(T)= \sum_{n:T_n\le T}T_n^{-1}\le C_0^{1/2}\sum_{n:T_n\le T}n^{-1/2}
\le C_0^{1/2} \sum_{n\le C_0T^2}n^{-1/2} \le
$$ $$
\le C_0^{1/2}\int_0^{C_0T^2}x^{-1/2}\,dx= C_0^{1/2}\cdot 2(C_0T^2)^{1/2}=
2C_0T.
$$

Set $T_0=l_m\sqrt{S}$, where $l_1=l_2=2$ and $l_m=2^{2^{4m}}$ for $m\ge3$.
Let $\La$ be a periodic cylinder and $|\La|$ be the length of $\La$.  For
any $\la\ge1$ let $A_\La(\la)$ denote the set of directions $v\in S^1$ such
that the projection of periodic geodesics from $\La$ on the direction
orthogonal to $v$ is at most $\la^{-1}T_0$.  Let $v\in A_\La(\la)$ and
$\varphi$ be the angle between $v$ and the direction of the cylinder $\La$
($0\le\varphi\le\pi/2$).  Then $\varphi\le \pi/2\cdot\sin\varphi\le
\pi/2\cdot\la^{-1}T_0/|\La|$.  It follows that
$$
\nu(A_\La(\la))\le 4\cdot \pi/2\cdot\la^{-1}T_0/|\La|=
2\pi T_0\la^{-1}\cdot |\La|^{-1},
$$
where $\nu$ is Lebesgue measure on the circle $S^1$ normalized so that
$\nu(S^1)=2\pi$.

Take an arbitrary number $T\ge T_0$ and set $\la=T/T_0$.  Let $\om$ denote
the translation structure of the translation surface $M$.  For any
direction $v\in S^1$ define an operator $a_{\la,v}\in\SL$ by equalities
$a_{\la,v}v=\la^{-1}v$, $a_{\la,v}u=\la u$, where $u\in S^1$ is a vector
orthogonal to $v$.  We claim that there exists a periodic cylinder $\La$ of
area at least $S/m$ such that the length of $\La$ with respect to the
translation structure $a_{\la,v}\,\om$ does not exceed $T_0$.  In the case
$m\ge3$, this follows from Theorem \ref{main2}.  In the case $m\le2$, the
translation surface $M$ is a torus, thus the claim follows from Lemma
\ref{lower3}.  The projection of a geodesic from $\La$ on the direction $u$
(with respect to the translation structure $\om$) is at most $\la^{-1}T_0$,
hence $v\in A_\La(\la)$.  Since $\la\ge1$, the condition number (see
Section \ref{exist}) of the operator $a_{\la,v}$ is equal to $\la$.
Therefore, $|\La|\le\la T_0=T$.  It follows that the union of sets
$A_\La(\la)$ over all periodic cylinders $\La$ of length at most $T$ (with
respect to $\om$) and of area at least $S/m$ is the circle $S^1$.

Let $\Sigma$ denote the sum of measures of sets $A_\La(\la)$, where $\La$
runs over all periodic cylinders of length at most $T$ and of area at least
$S/m$.  Since these sets cover the whole circle $S^1$, we have $\Sigma\ge
2\pi$.  Set $\alpha=(4T_0^2C_0)^{-1}$.  By Theorem \ref{main2} and Lemma
\ref{lower3}, the translation structure $\om$ admits a periodic cylinder of
length at most $T_0$, hence $1\le N_1(M,T_0)\le C_0T_0^2$.  In particular,
$\alpha\le1/4$.  The sum $\Sigma$ can be written as $\Sigma_1+\Sigma_2$,
where $\Sigma_1$ is the sum of summands corresponding to the cylinders of
length at most $\alpha T$, and $\Sigma_2$ is the sum over cylinders of
length greater than $\alpha T$.  It follows from the above estimate of
$\nu(A_\La(\la))$ that
$$
\Sigma_1\le 2\pi T_0\la^{-1}\cdot \sigma(\alpha T)\le
2\pi T_0\la^{-1}\cdot 2C_0\alpha T, \qquad
\Sigma_2\le 2\pi T_0\la^{-1}\cdot (\alpha T)^{-1} N(T),
$$
where $N(T)$ is the number of periodic cylinders of length at most $T$ and
of area at least $S/m$.  Then
$$
2\pi\le 2\pi T_0\la^{-1} (2C_0\alpha T+ (\alpha T)^{-1} N(T)),
$$
thus,
$$
N(T)\ge (T_0^{-1}\la- 2C_0\alpha T)(\alpha T)= (8T_0^4C_0)^{-1}T^2.
$$
Consequently,
$$
N_2(M,T)/S\ge N(T)/m\ge (8mT_0^4C_0)^{-1}T^2=
(8ml_m^4\tilde h_m)^{-1}s^2S^{-2}T^2.
$$

To complete the proof, it remains to show that $8ml_m^4\tilde h_m\le
(600m)^{(2m)^{2m}}$.  If $m=1$, then $8ml_m^4\tilde h_m=7\cdot 2^7<600^4$.
If $m=2$, then $8ml_m^4\tilde h_m=2^8\cdot 800^{4^4}<1200^{4^4}$.  In the
case $m\ge3$, we have $(m/2)^{2m}\ge (3/2)^6>10$.  Then $(2m)^{2m}>10\cdot
4^{2m}=10\cdot 2^{4m}$.  It follows that $(3/2)^{(2m)^{2m}}>2^{5\cdot
2^{4m}}$.  Besides, $2^{4m}>8m$.  Finally,
$$
(8ml_m^4\tilde h_m)^{-1} (600m)^{(2m)^{2m}}=
(8m)^{-1}\cdot 2^{-4\cdot 2^{4m}}\cdot (3/2)^{(2m)^{2m}}>
$$ $$
(8m)^{-1}\cdot 2^{2^{4m}}>2^{2^{4m}}\cdot 2^{-4m}>2^{4m}>1.
$$
The theorem is proved.
\end{proofof}

\section{Moduli spaces of translation structures}\label{quad}

In this section we consider moduli spaces of translation structures on a
given surface and properties of periodic geodesics of a generic translation
structure.

Let $M$, $M'$ be translation surfaces, and $\om$, $\om'$ be their
translation structures.  An orientation-preserving homeomorphism $f:M\to
M'$ is called an {\em isomorphism\/} of the translation surfaces if $f$
maps the set of singular points of $M$ onto the set of singular points of
$M'$ and $f$ is a translation in local coordinates of the atlases $\om$ and
$\om'$.  The translation structures $\om$ and $\om'$ are called {\em
isomorphic\/} if there is an isomorphism $f:M\to M'$.  If the isomorphism
can be chosen isotopic to the identity, then the structures $\om$ and
$\om'$ are called {\em isotopic}.  A homeomorphism $f:M\to M'$ is said to
be {\em piecewise affine\/} if there exists a triangulation of the
translation surface $M$ such that $f$ is affine on every triangle of the
triangulation in local coordinates of the atlases $\om$ and $\om'$.  The
linear parts $a_1,\dots,a_k$ of restrictions of $f$ to the triangles are
uniquely determined.  Set $b(f)= \max(\|a_1-1\|,\dots,\|a_k-1\|)$.
Clearly, $b(f)=0$ if and only if $f$ is an isomorphism of translation
surfaces.

Given positive integers $p$ and $n$, let $M_p$ be a compact connected
oriented surface of genus $p$ and $Z_n$ be a subset of $M_p$ of cardinality
$n$.  Denote by $\Omega(p,n)$ the set of translation structures on $M_p$
such that $Z_n$ is the set of singular points.  By $\tcM(p,n)$ denote the
set of equivalence classes of isotopic translation structures in
$\Omega(p,n)$, and by $\cM(p,n)$ denote the set of equivalence classes of
isomorphic translation structures in $\Omega(p,n)$.  Both sets $\tcM(p,n)$
and $\cM(p,n)$ can serve as {\em moduli spaces\/} of translation structures
on $M_p$ with $n$ singular points.

Given a translation structure $\om=\{(U_\alpha,\phi_\alpha)\}_{\alpha\in
\cA}$ on $M_p$ and a homeomorphism $f:M_p\to M_p$, the atlas $\om f=
\{(f^{-1}(U_\alpha),\phi_\alpha f)\}_{\alpha\in\cA}$ is a translation
structure on $M_p$ isomorphic to $\om$.  Let $H(p,n)$ denote the group of
orientation-preserving homeomorphisms of the surface $M_p$ leaving
invariant the set $Z_n$.  By $H_0(p,n)$ denote the subgroup of $H(p,n)$
consisting of homeomorphisms isotopic to the identity.  For any
$\om\in\Omega(p,n)$ and $f\in H(p,n)$, the translation structure $\om f$
belongs to $\Omega(p,n)$.  The map $H(p,n)\times\Omega(p,n)\ni(f,\om)
\mapsto \om f^{-1}$ defines an action of the group $H(p,n)$ on the set
$\Omega(p,n)$.  By definition, $\tcM(p,n)=\Omega(p,n)/H_0(p,n)$ and
$\cM(p,n)=\Omega(p,n)/H(p,n)$.  The modular group $\Mod(p,n)=
H(p,n)/H_0(p,n)$ acts naturally on the set $\tcM(p,n)$ and $\cM(p,n)=
\tcM(p,n)/\Mod(p,n)$.  Further, the group $\SL$ acts on the set
$\Omega(p,n)$ as defined in Section \ref{prel}.  Obviously, this action
commutes with the action of $H(p,n)$.  Therefore the action of $\SL$
descends to actions on the spaces $\tcM(p,n)$ and $\cM(p,n)$.  The action
of the group $\SL$ on $\tcM(p,n)$ commutes with the action of $\Mod(p,n)$.

For any $\om\in\Omega(p,n)$, let $M_p(\om)$ denote the surface $M_p$
considered as the translation surface with the translation structure $\om$.
Suppose $\om\in\Omega(p,n)$ and $\eps>0$.  By definition, a translation
structure $\om'\in\Omega(p,n)$ belongs to the set $W(\om,\eps)$ if there
exists a piecewise affine map $f:M_p(\om)\to M_p(\om')$ such that
$b(f)<\eps$ and $f\in H(p,n)$.  Further, $\om'\in\widetilde W(\om,\eps)$ if
the map $f$ can be chosen in $H_0(p,n)$.  The collection of sets
$W(\om,\eps)$, where $\om\in\Omega(p,n)$ and $\eps>0$, generates a topology
$\cW$ on $\Omega(p,n)$. The collection of sets $\widetilde W(\om,\eps)$
generates a stronger topology $\tcW$.  The topology $\tcW$ descends to a
Hausdorff topology on $\tcM(p,n)$, while the topology $\cW$ descends to a
Hausdorff topology on $\cM(p,n)$.  It can be shown that the topologies on
$\tcM(p,n)$ and $\cM(p,n)$ are topologies of locally compact metric spaces.
The group $\SL$ acts on the spaces $\tcM(p,n)$ and $\cM(p,n)$ by
homeomorphisms.  The action of $\Mod(p,n)$ on $\tcM(p,n)$ is also by
homeomorphisms, besides, this action is properly discontinuous.

Let $\gamma:[0,1]\to M_p$ be a continuous curve.  For any $\om\in
\Omega(p,n)$ there exists a continuous curve $\gamma_\om:[0,1]\to\bR^2$
such that $\gamma$ is a translation of $\gamma_\om$ in coordinates of the
atlas $\om$.  The curve $\gamma_\om$ is determined up to translation.  The
vector $\gamma_\om(1)-\gamma_\om(0)$ is called the {\em holonomy\/} of the
curve $\gamma$ with respect to translation structure $\om$ and is denoted
by $\hol_\om(\gamma)$.  If $\gamma$ is a geodesic segment of $\om$, then
the vector $\hol_\om(\gamma)$ is of the same length and direction as
$\gamma$.  The holonomy $\hol_\om(\gamma)$ does not change if we replace
the curve $\gamma$ by a homologous one or replace the translation structure
$\om$ by an isotopic one.  In particular, the map $\hol_\om$ is
well-defined for $\om\in\tcM(p,n)$.  Also, the map $\hol_\om$ extends to a
map of the relative homology group $H_1(M_p,Z_n;\bZ)$ that is an element of
the relative cohomology group $H^1(M_p,Z_n;\bR^2)$.

Suppose $\Gamma=(\gamma_1,\dots,\gamma_N)$ is an ordered basis of the group
$H_1(M_p,Z_n;\bZ)$.  Note that $N=2p+n-1$.  Define a map
$C_\Gamma:\tcM(p,n)\to(\bR^2)^N\approx\bR^{2N}$ by $C_\Gamma(\om)=
(\hol_\om(\gamma_1),\dots,\hol_\om(\gamma_N))$.  The map $C_\Gamma$ is a
local homeomorphism (see \cite{Veech90}).  For any ordered basis $\Gamma'$
of $H_1(M_p,Z_n;\bZ)$ there exists a unique linear operator $g\in
\GL(2N,\bR)$ such that $C_{\Gamma'}=gC_\Gamma$.  It is easy to observe that
$g\in\GL(2N,\bZ)$. The inverse operator $g^{-1}$ is also in $\GL(2N,\bZ)$,
hence $|\det g|=1$. Thus the collection of maps of the form $C_\Gamma$
endows the space $\tcM(p,n)$ with the structure of a real analytic
$2N$-dimensional manifold along with a volume element.  Every element
$\phi\in\Mod(p,n)$ induces an automorphism $\phi_*$ of the group
$H_1(M_p,Z_n;\bZ)$.  Clearly, $C_\Gamma(\om\phi)=C_{\Gamma'}(\om)$ for any
$\om\in\tcM(p,n)$, where $\Gamma'=(\phi_*^{-1}\gamma_1,\dots,
\phi_*^{-1}\gamma_N)$.  It follows that the action of $\Mod(p,n)$ on the
space $\tcM(p,n)$ is analytic and preserves the volume element.  Further,
$C_\Gamma(g\om)=(g\hol_\om(\gamma_1),\dots,g\hol_\om(\gamma_N))$ for any
$g\in\SL$ and $\om\in\tcM(p,n)$, hence the action of $\SL$ on $\tcM(p,n)$
is also real analytic and also preserves the volume element.

For any $\om\in\tcM(p,n)$, let $a(\om)$ denote the area of the surface
$M_p$ with respect to translation structures in the isotopy class $\om$. It
is easy to see that $a(\om)$ is a quadratic form of the vector
$C_\Gamma(\om)$.  Therefore the set $\tcM_1(p,n)=a^{-1}(1)$ is a real
analytic submanifold of $\tcM(p,n)$ of codimension $1$.  This submanifold
is invariant under the actions of $\Mod(p,n)$ and $\SL$.  The volume
element on $\tcM(p,n)$ induces a volume element on $\tcM_1(p,n)$.  By
$\tilde\mu_0$ denote the corresponding Borel measure on $\tcM_1(p,n)$.  Let
$\pi_0:\tcM(p,n)\to\cM(p,n)$ be the natural projection.  The set
$\cM_1(p,n)=\pi_0(\tcM_1(p,n))$ is a topological subspace of $\cM(p,n)$
invariant under the action of $\SL$.  It can be shown that the number of
connected components of $\cM_1(p,n)$ is at most finite.  Since the action
of the group $\Mod(p,n)$ on $\tcM_1(p,n)$ is properly discontinuous, there
exists a unique Borel measure $\mu_0$ on $\cM_1(p,n)$ such that
$\mu_0(\pi_0(U))=\tilde\mu_0(U)$ whenever the set $U\subset\tcM_1(p,n)$ is
Borel and $\pi_0$ is injective on $U$.  The measure $\mu_0$ is invariant
under the action of $\SL$.  In addition, the measure $\mu_0$ is finite (see
\cite{Masur82}, \cite{Veech90}).

The maps $H(p,n)\times\Omega(p,n)\times M_p\ni(f,\om,x)\mapsto(\om
f^{-1},f(x))$ and $\SL\times\Omega(p,n)\times M_p\ni(g,\om,x)\mapsto
(g\om,x)$ define commuting actions of the groups $H(p,n)$ and $\SL$ on the
set $\Omega(p,n)\times M_p$.  Set $\tcY(p,n)=(\Omega(p,n)\times
M_p)/H_0(p,n)$ and $\cY(p,n)=(\Omega(p,n)\times M_p)/H(p,n)$.  The
topologies $\tcW$ and $\cW$ on the set $\Omega(p,n)$ give rise to
topologies $\tcW_1$ and $\cW_1$ on $\Omega(p,n)\times M_p$.  The topologies
$\tcW_1$ and $\cW_1$ descend to Hausdorff topologies on $\tcY(p,n)$ and
$\cY(p,n)$, respectively.  Let $\tilde p_0:\tcY(p,n)\to\tcM(p,n)$ and
$p_0:\cY(p,n)\to\cM(p,n)$ be the natural projections.  The group $\SL$ acts
on the spaces $\tcY(p,n)$ and $\cY(p,n)$ by homeomorphisms.  The subspaces
$\tcY_1(p,n)=\tilde p_0^{-1}(\tcM_1(p,n))$ and $\cY_1(p,n)=
p_0^{-1}(\cM_1(p,n))$ are invariant under these actions.

Let $\tilde\om\in\tcM(p,n)$.  Take a translation structure $\om$ in the
isotopy class $\tilde\om$.  Fix a triangulation $\tau$ of the surface $M_p$
by pairwise disjoint saddle connections of $\om$.  By definition, a
translation structure $\om'\in\Omega(p,n)$ belongs to the set $X(\om,\tau)$
if there exists a piecewise affine map $f:M_p(\om)\to M_p(\om')$ such that
$f\in H_0(p,n)$ and $f$ is affine on every triangle of $\tau$.  Let
$\widetilde X(\om,\tau)$ denote the set of isotopy classes $\tilde\om'\in
\tcM(p,n)$ that have representatives in $X(\om,\tau)$.  The set $\widetilde
X(\om,\tau)$ is open and each $\tilde\om'\in\widetilde X(\om,\tau)$ has
precisely one representative in $X(\om,\tau)$.  Hence each $\eta\in\tilde
p_0^{-1}(\widetilde X(\om,\tau))$ has precisely one representative in
$X(\om,\tau)\times M_p$.  This gives rise to a map $F_{\om,\tau}:\tilde
p_0^{-1}(\widetilde X(\om,\tau))\to M_p$, which is continuous.  The map
$\tilde p_0^{-1}(\widetilde X(\om,\tau))\ni\eta\mapsto(\tilde p_0(\eta),
F_{\om,\tau}(\eta))$ is a homeomorphism of $\tilde p_0^{-1}(\widetilde
X(\om,\tau))$ onto $\widetilde X(\om,\tau)\times M_p$.  It follows that the
space $\tcY(p,n)$ is a fiber bundle over $\tcM(p,n)$ with the fiber $M_p$.

For any $\om\in\cM(p,n)$ the following conditions are equivalent: (i)
translation structures in the isomorphy class $\om$ have no automorphisms
different from the identity; (ii) for any $\tilde\om\in\pi_0^{-1}(\om)$ the
restriction of the projection $\pi_0$ to some neighborhood of $\tilde\om$
is a homeomorphism.  Let $U_0$ be the open set of $\om\in\cM(p,n)$
satisfying these conditions.  The preimage $p_0^{-1}(U_0)\subset\cY(p,n)$
is a fiber bundle over $U_0$ with the fiber $M_p$.  Suppose $\om\in
\cM(p,n)\setminus U_0$ and $\om_0\in\om$; then $p_0^{-1}(\om)$ is
homeomorphic to $M_p/\Aut(\om_0)$, where $\Aut(\om_0)$ is the group of
automorphisms of the translation structure $\om_0$.  Since $U_0$ is an open
dense subset of full measure of $\cM(p,n)$, we consider $\cY(p,n)$ as a
fiber bundle over $\cM(p,n)$ with the fiber $M_p$ (even though some fibers
may be not homeomorphic to $M_p$).

For any $\om\in\Omega(p,n)$, let $\xi_\om$ denote the Borel measure on
$M_p$ induced by the translation structure $\om$.  Let $\tilde\om$ be the
isotopy class of $\om$.  The map $h_\om:M_p\to\tilde p_0^{-1}(\tilde\om)$
defined by the relation $(\om,x)\in h_\om(x)$, $x\in M_p$, is a
homeomorphism.  The measure $\tilde\nu_{\tilde\om}=\xi_\om h_\om^{-1}$ on
the fiber $\tilde p_0^{-1}(\tilde\om)$ does not depend on the choice of
$\om\in\tilde\om$.  Likewise, for any $\tilde\om\in\cM(p,n)$ the measures
on $M_p$ induced by translation structures in the isomorphy class $\tilde
\om$ define a Borel measure $\nu_{\tilde\om}$ on $p_0^{-1}(\tilde\om)$
(even if the fiber $p_0^{-1}(\tilde\om)$ is not homeomorphic to $M_p$).
The space $\tcY_1(p,n)$, which is a fiber bundle over $\tcM_1(p,n)$ with
the fiber $M_p$, carries a natural measure $\tilde\mu_1$ that is the
measure $\tilde\mu_0$ on the base $\tcM_1(p,n)$ and is the measure
$\tilde\nu_{\tilde\om}$ on the fiber $\tilde p_0^{-1}(\tilde\om)$.  In
other words, $d\tilde\mu_1(\eta)=d\tilde\nu_{\tilde\om}(\eta)\,d\tilde
\mu_0(\tilde\om)$.  Similarly, the space $\cY_1(p,n)$ carries a natural
measure $\mu_1$ such that $d\mu_1(\eta)=d\nu_{\tilde\om}(\eta)\,
d\mu_0(\tilde\om)$.  The measures $\tilde\mu_1$ and $\mu_1$ are invariant
under the actions of the group $\SL$ on the spaces $\tcY_1(p,n)$ and
$\cY_1(p,n)$, respectively.  Let $\pi_1:\tcY(p,n)\to\cY(p,n)$ be the
natural projection.  Then $\mu_1(\pi_1(U))=\tilde\mu_1(U)$ for any Borel
set $U\subset\tcY_1(p,n)$ such that $\pi_1$ is injective on $U$.

We proceed to the proof of Theorems \ref{main6} and \ref{main7}.  In what
follows constructions and results of the papers of Veech \cite{Veech98} and
of Eskin and Masur \cite{EM} are used extensively.  The proof relies on
Theorem \ref{quad1} formulated below.  The formulation requires additional
definitions and notation.

Let $V$ denote a pair of sequences $v_1,v_2,\dots$ and $w_1,w_2,\dots$,
where elements of the first sequence are nonzero vectors in $\bR^2$ and
elements of the second sequence are positive numbers.  The number $w_k$ is
called the {\em weight\/} of the vector $v_k$.  It is assumed that the
sequence of vectors tends to infinity or is finite, and the sequence of
weights is bounded.  By $\cV$ denote the set of all such pairs.  Two pairs
$V_1,V_2\in\cV$ are considered to be equal if one of them can be obtained
from the other by reordering its vectors along with the corresponding
reorder of weights.  The group $\SL$ acts on the set $\cV$ by the natural
action on vectors and the trivial action on weights.  To each pair
$V\in\cV$ we assign a linear functional $\Phi[V]$ on the space $C_c(\bR^2)$
of continuous compactly supported functions on $\bR^2$; the functional is
defined by the relation $\Phi[V](f)=\sum_{k=1}^\infty w_kf(v_k)$.  Note
that two elements $V_1,V_2\in\cV$ are equal if and only if $\Phi[V_1]=
\Phi[V_2]$.  Furthermore, for any $T>0$ set $N_V(T)=\sum_{k:|v_k|\le T}
w_k$.  The function $N_V$ is called the {\em growth function\/} of $V$.

Let $\cM$ be a locally compact metric space endowed with a finite Borel
measure $\mu$.  Suppose the group $\SL$ acts on the space $\cM$ by
homeomorphisms.  We assume that the measure $\mu$ is invariant under this
action and the action is {\em ergodic}, that is, any measurable subset of
$\cM$ invariant under the action is of zero or full measure.  Let $V$ be a
map of the space $\cM$ to $\cV$.  The map $V$ is supposed to satisfy the
following conditions:

(0) for any $f\in C_c(\bR^2)$ the function $\cM\ni\om\mapsto
\Phi[V(\om)](f)$ is Borel;

(A) the map $V$ intertwines the actions of the group $\SL$ on the spaces
$\cM$ and $\cV$, that is, $V(g\om)=gV(\om)$ for any $g\in\SL$ and any
$\om\in\cM$;

(B) for any $\om\in\cM$ there exists a constant $c=c(\om)>0$ such that
$N_{V(\om)}(T)\le cT^2$ for $T>1$; the constant $c$ can be chosen uniformly
as $\om$ varies over a compact subset of $\cM$;

(C) there exist positive constants $T_0$ and $\eps$ such that the function
$\om\mapsto N_{V(\om)}(T_0)$ belongs to the space $L^{1+\eps}(\cM,\mu)$.

\begin{theorem}\label{quad1}
Suppose a map $\cM\ni\om\mapsto V(\om)\in\cV$ satisfies conditions (0),
(A), (B), and (C).  Then (a) for any $f\in C_c(\bR^2)$ the function
$\cM\ni\om\mapsto \Phi[V(\om)](f)$ is integrable and
$$
\frac1{\mu(\cM)}\int_{\cM}\Phi[V(\om)](f)\,d\mu(\om)=
c_V\int_{\bR^2}f(x)\,dx,
$$
where $c_V$ is a nonnegative constant;
\par (b) for $\mu$-almost every $\om\in\cM$ one has
$$
\lim_{T\to\infty} N_{V(\om)}(T)/T^2=\pi c_V.
$$
\end{theorem}

The first statement of Theorem \ref{quad1} was proved by Veech
\cite{Veech98}.  He also proved that
$$
\lim_{T\to\infty}\int_{\cM}\, \Bigl| N_{V(\om)}(T)/T^2-\pi c_V \Bigr|
\,d\mu(\om)=0.
$$
The second statement was proved by Eskin and Masur \cite{EM}.

\begin{proofof}{Theorems \ref{main6} and \ref{main7}}
Let $\cC$ be a connected component of the space $\cM_1(p,n)$ ($p,n\ge1$)
and $Y$ be the component of $\cY_1(p,n)$ that is a fiber bundle over
$\cC$ with respect to the natural projection $p_0:\cY_1(p,n)\to\cM_1(p,n)$.
In what follows we often regard elements of $\cM(p,n)$ and $\tcM(p,n)$ as
translation structures on $M_p$, and elements of $\cY(p,n)$ and $\tcY(p,n)$
as pairs in $\Omega(p,n)\times M_p$ (although, in fact, all such elements
are equivalence classes).  We define maps $V_1:\cC\to\cV$, $V_2:\cC\to\cV$,
and $V_3:Y\to\cV$ as follows.  For any $\om\in\cC$ the pairs $V_1(\om)$ and
$V_2(\om)$ share the same sequence of vectors that is the sequence of
vectors associated to periodic cylinders of the translation structure
$\om$.  Note that to any periodic cylinder we associate two vectors of the
same length and of opposite directions.  Both vectors are supposed to be in
the sequence.  If a vector is associated to $k>1$ different periodic
cylinders, it is to appear $k$ times in the sequence.  All weights of
$V_1(\om)$ are equal to $1$.  For $V_2(\om)$, the weight of a vector
associated to a periodic cylinder is the area of the cylinder.  Further,
for any $x\in M_p$ the sequence of vectors of $V_3(\om,x)$ is defined to be
the sequence of vectors associated to periodic geodesics of $\om$ passing
through $x$.  All weights of $V_3(\om,x)$ are equal to $1$.  By definition,
$N_{V_1(\om)}(T)=2N_1(\om,T)$, $N_{V_2(\om)}(T)=2N_2(\om,T)$, and
$N_{V_3(\om,x)}(T)=2N_3(\om,x,T)$ for any $T>0$.

Let us show that the maps $V_1$, $V_2$, and $V_3$ satisfy all hypotheses of
Theorem \ref{quad1}.  First notice that the natural actions of the group
$\SL$ on the spaces $\cC$ and $Y$ are ergodic.  The ergodicity of the
action on $\cC$ was proved by Veech \cite{Veech86}, and the ergodicity of
the action on $Y$ was proved by Eskin and Masur \cite{EM}.  By definition
of the actions of $\SL$ on the spaces $\cC$, $Y$, and $\cV$, the maps
$V_1$, $V_2$, and $V_3$ satisfy condition (A).

Let $\cS(p,n)$ be the set of free holonomy classes of simple closed
oriented curves in $M_p\setminus Z_n$.  For any $\gamma\in\cS(p,n)$, the
map $\tcM(p,n)\ni\om\mapsto\hol_\om(\gamma)\in\bR^2$ is continuous.  By
$U(\gamma)$ denote the set of translation structures in $\tcM(p,n)$ that
admit a periodic geodesic in the holonomy class $\gamma$.  By $U_1(\gamma)$
denote the set of pairs $(\om,x)\in\tcY(p,n)$ such that some periodic
geodesic of the translation structure $\om$ passing through the point $x$
is in the holonomy class $\gamma$.  Obviously, the sets $U(\gamma)$ and
$U_1(\gamma)$ are open.  Given $\om\in U(\gamma)$, all periodic geodesics
of $\om$ that belong to the holonomy class $\gamma$ form one periodic
cylinder.  Let $a_\gamma(\om)$ denote the area of this cylinder.  For any
$\om\notin U(\gamma)$, put $a_\gamma(\om)=0$.  It is easy to observe that
the function $a_\gamma$ is continuous on $\tcM(p,n)$.  Let $\pi_0:
\tcM(p,n)\to\cM(p,n)$ and $\pi_1:\tcY(p,n)\to\cY(p,n)$ be the canonical
projections.  Then for any $\om\in\pi_0^{-1}(\cC)$, any
$\eta\in\pi_1^{-1}(Y)$, and any $f\in C_c(\bR^2)$ we have
\begin{eqnarray*}
\Phi[V_1(\pi_0(\om))](f) &=& \sum_{\gamma\in\cS(p,m)}
\chi_{U(\gamma)}(\om)\, f(\hol_\om(\gamma)),\\
\Phi[V_2(\pi_0(\om))](f) &=& \sum_{\gamma\in\cS(p,m)}
a_\gamma(\om)\, f(\hol_\om(\gamma)),\\
\Phi[V_3(\pi_1(\eta))](f) &=& \sum_{\gamma\in\cS(p,m)}
\chi_{U_1(\gamma)}(\eta)\, f(\hol_{\tilde p_0(\eta)}(\gamma)).
\end{eqnarray*}
All three sums are locally finite.  It follows that the functions
$\pi_0^{-1}(\cC)\ni\om\mapsto\Phi[V_1(\pi_0(\om))](f)$ and
$\pi_1^{-1}(Y)\ni\eta\mapsto\Phi[V_3(\pi_1(\eta))](f)$ are Borel, while the
function $\pi_0^{-1}(\cC)\ni\om\mapsto\Phi[V_2(\pi_0(\om))](f)$ is
continuous.  Then the functions $\cC\ni\om\mapsto\Phi[V_1(\om)](f)$ and
$Y\ni\eta\mapsto\Phi[V_3(\eta)](f)$ are also Borel and the function
$\cC\ni\om\mapsto\Phi[V_2(\om)](f)$ is also continuous.  Thus condition (0)
holds for the maps $V_1$, $V_2$, and $V_3$.

For any $\om\in\cC$, let $s(\om)$ denote the length of the shortest saddle
connection of the translation structure $\om$.  The function $\om\mapsto
s(\om)$ is continuous and bounded on $\cC$.  Therefore the upper estimate
in Theorem \ref{main5} implies the map $V_1$ satisfies condition (B).  To
verify condition (C), we need the following theorem.

\begin{theorem}[\cite{EM}]\label{quad2}
(1) Given $T>0$ and $\eps>0$, there exists a positive constant $C_{T,\eps}$
such that
$$
N_1(\om,T)\le C_{T,\eps}(s(\om))^{-1-\eps}
$$
for any $\om\in\cC$.

(2) For any $\beta\in[1,2)$ the function $s^{-\beta}$ belongs to the space
$L^1(\cC,\mu_0)$.
\end{theorem}

Theorem \ref{quad2} implies that condition (C) holds for the map $V_1$.
Let $\om$ be an arbitrary translation structure in $\cC$.  By definition,
$N_2(\om,T)\le N_1(\om,T)$ for any $T>0$, and $N_3(\om,x,T)\le N_1(\om,T)$
for any $x\in M_p$ and any $T>0$.  It follows that conditions (B) and (C)
are satisfied by the maps $V_2$ and $V_3$ whenever these conditions are
satisfied by $V_1$.

Now it follows from Theorem \ref{quad1} that there exist constants
$c_1(\cC),c_2(\cC),c_3(\cC)\ge0$ such that $\lim\limits_{T\to\infty}
N_1(\om,T)/T^2=c_1(\cC)$ and $\lim\limits_{T\to\infty}N_2(\om,T)/T^2=
c_2(\cC)$ for $\mu_0$-almost every $\om\in\cC$, and
$\lim\limits_{T\to\infty}N_3(\om,x,T)/T^2=c_3(\cC)$ for $\mu_1$-almost
every $(\om,x)\in Y$.  The positivity of the numbers $c_1(\cC)$ and
$c_2(\cC)$ follows from the lower estimates in Theorem \ref{main5}.  It
remains to prove that $c_3(\cC)=c_2(\cC)$.  Take a function $f_0\in
C_c(\bR^2)$ such that $\int f_0(x)\,dx=1$.  By Theorem \ref{quad1}, we have
\begin{eqnarray*}
\frac1{\mu_0(\cC)}\int_{\cC}\Phi[V_2(\om)](f_0)\,d\mu_0(\om) &=&
2\pi^{-1}c_2(\cC),\\
\frac1{\mu_1(Y)}\int_Y\Phi[V_3(\eta)](f_0)\,d\mu_1(\eta) &=&
2\pi^{-1}c_3(\cC).
\end{eqnarray*}
For any $\om\in\cC$, let $\nu_\om$ denote the Borel measure on the fiber
$p_0^{-1}(\om)$ induced by translation structures in the equivalence class
$\om$.  It is easy to observe that
$$
\Phi[V_2(\om)](f)=\int_{p_0^{-1}(\om)} \Phi[V_3(\eta)](f)\,d\nu_\om(\eta)
$$
for any $f\in C_c(\bR^2)$.  Then
$$
\int_Y \Phi[V_3(\eta)](f_0)\,d\mu_1(\eta)= \int_{\cC}
\int_{p_0^{-1}(\om)}\Phi[V_3(\eta)](f_0)\,d\nu_\om(\eta)\,d\mu_0(\om)=
\int_{\cC}\Phi[V_2(\om)](f_0)\,d\mu_0(\om),
$$
besides,
$$
\mu_1(Y)=\int_{\cC}\nu_\om(p_0^{-1}(\om))\,d\mu_0(\om)=\mu_0(\cC).
$$
Hence, $c_3(\cC)=c_2(\cC)$.  The theorems are proved.
\end{proofof}

\end{document}